\documentclass[11pt]{article}
\usepackage{latexsym,amssymb}
\def\demo{\noindent{\bf Proof. }}
\def\QED{\hfill$\Box$}
\newtheorem{Theorem}{Theorem}[section]
\newtheorem{Lemma}[Theorem]{Lemma}
\newtheorem{Corollary}[Theorem]{Corollary}
\newtheorem{Proposition}[Theorem]{Proposition}

\newtheorem{Definition}[Theorem]{Definition}

\begin{document}
\topmargin3mm
\hoffset=-1cm
\voffset=-1.5cm
\

\medskip
\begin{center}
{\large\bf On the Stanley depth of edge ideals of k--partite clutters
}
\vspace{6mm}\\
\footnotetext{
\\
Partially supported by CONACYT}
\end{center}

\medskip
\begin{center}
Luis A. Dupont and Daniel G. Mendoza
\\
{\small Facultad de Matem\'aticas, Universidad Veracruzana}\vspace{-1mm}\\
{\small Circuito Gonzalo Aguirre Beltr\'an S/N;}\vspace{-1mm}\\
{\small Zona Universitaria;}\vspace{-1mm}\\
{\small Xalapa, Ver., M\'exico, CP 91090.}\vspace{-1mm}\\
{\small e-mail: {\tt ldupont@uv.mx; dgmendozaramirez@gmail.com}\vspace{4mm}}
\end{center}

\medskip

\begin{abstract}
\noindent
We give upper bounds for the Stanley depth of edge ideals of certain $k$--partite clutters.
In particular, we genera\-li\-ze a result of Ishaq about the Stanley depth of the edge ideal
of a complete bipartite graph.
A result of Pournaki, Seyed Fakhari and Yassemi  implies that the Stanley's conjecture 
holds for d-uniform complete d-partite clutters. Here we give a shorter and different proof  of this fact.
\end{abstract}

\section{Introduction}
Let $R=K[x_1,\ldots,x_n]$ be a polynomial ring over a field $K.$
A {\it clutter\/} $\mathcal{C}$, with finite vertex set
$V=\{x_{1},...,x_{n}\}$ is a family of subsets of $V$, called edges, none
of which is included in another. The set of
vertices and edges of $\mathcal{C}$ are denoted by $V(\mathcal{C})$
and $E(\mathcal{C})$ respectively. For example, a simple
graph (no multiple edges or loops) is a clutter. The {\it edge ideal\/} of $\mathcal{C}$,
denoted by $I(\mathcal{C})$, is the ideal of $R$
generated by all monomials $x_e=\prod_{x_i\in e}x_i$ such
that $e\in E(\mathcal{C})$. The map
$$\mathcal{C}\longmapsto I(\mathcal{C})$$
gives a one to one
correspondence between the family of clutters and the family of
squarefree monomial ideals. Edge ideals of graphs were introduced and studied in \cite{ITG,Vi2}.
Edge ideals of clutters correspond to simplicial complexes via the
Stanley-Reisner correspondence \cite{Stanley} and to facet ideals
\cite{faridi,zheng}.
\\

A $k$--partite clutter is a clutter $\mathcal{C}$  where the vertices are partitioned into k subsets $V(\mathcal{C}) = V_{1} \cup V_{2} \cup \cdots \cup V_{k}$ with the following conditions:
\\
$(1)$ No two vertices in the same subset are adjacent, i.e.,  $|V_{i} \cap E|\leq 1$ for all $1\leq i \leq k$ and $E \in E(\mathcal{C})$.
\\
$(2)$ There is no partition of the vertices with fewer than k subsets where condition $(1)$ holds.
\\
A clutter is called $d$--{\it uniform\/} or {\it uniform\/} if all its edges have
exactly $d$ vertices. Along the paper we introduce most of the notions that are relevant for
our purposes. Our main re\-fe\-rences for combinatorial optimization and commutative
algebra are \cite{cornu-book,tesis-docto-dupont,monalg}.
\\

Let $M$ be a finitely generated $\mathbb{Z}^{n}$-graded $R$-module, $R=K[x_1,\ldots,x_n]$.
If $u\in M$ is a homogeneous element in $M$ and $Z\subseteq\{x_{1}, ..., x_{n}\}$ then let $uK[Z]\subset M$ denote the linear $K$--subspace of $M$ of all elements $uf$, $f\in K[Z]$. This space is called a {\it Stanley space\/} of
dimension $|Z|$ if $uK[Z]$ is a free $K[Z]$-module.
A presentation of $M$ as a finite direct sum of Stanley spaces
$$\mathcal{D} : M = \bigoplus_{i=1}^{r} u_{i}K[Z_{i}]$$
is called a {\it Stanley decomposition\/} of $M$. The number
$$ sdepth(\mathcal{D}) = \textmd{min}\{ |Z_{i}| : i=1,\ldots , r\}$$
is called the Stanley depth of decomposition $\mathcal{D}$ and the number
$$sdepth(M) := \textmd{max}\{sdepth(\mathcal{D}) : \mathcal{D} \textmd{  is a Stanley decomposition of  } M\}\leq n.$$
is called the {\it Stanley depth\/} of $M$. This is a combinatorial invariant which does
not depend on the characteristic of $K$.
\\

In $1982$, \cite{conj-Stanley}, Stanley introduced the idea of what is now called the {\it Stanley depth\/} of a $\mathbb{Z}^{n}$-graded module over a commutative ring and conjectured that $sdepth(M)\geq depth(M)$. While some special cases of the conjecture have been resolved, it still remains
largely open, (for example see \cite{Apel,Anwar-Popescu,Biro,Ishaq,Popescu1,Popescu2,Popescu3,Shen}). Shen's proof (see from \cite[Lema~2.3,Theorem~2.4]{Shen}) relies on a theorem of Cimpoea\c{s}, \cite[Theorem~2.1]{Cimpoeas}, which states that the Stanley depth of a complete intersection monomial ideal is equal to that of its radical, which allows for a focus on squarefree ideals. In \cite[Theorem~2.8]{Ishaq} Ishaq showed that the Stanley depth of the edge ideal of a complete bipartite graph over $n$ vertices with $n\geq4$ is less than or equal to $\frac{n+2}{2}$. In \cite{Ishaq2} Ishaq and Qureshi, provide an upper bound for the Stanley depth of an edge ideal of a $k$--uniform complete bipartite hypergraph which is a kind of generalization to the complete bipartite graph.
\\

The aim of this paper is to bound the Stanley depth of the edge ideal of a $d$--uniform complete $k$--partite clutter [Theorems~\ref{teorema1}, \ref{teorema1.seguido}, \ref{teorema2}]. The proofs use the co\-rres\-pondence between a Stanley decomposition of a monomial ideal and a partition of a particular poset into intervals established
by Herzog, Vladoiu, and Zheng. In \cite[Corollary 2.9]{pournaki} Pournaki, Seyed Fakhari and Yassemi show that the Stanley's conjecture holds for finite products of monomial prime ideals. This fact implies that the conjecture holds for d-uniform complete d-partite clutters. Here we give a shorter and different proof  of this result [Theorem~\ref{conjetura-stanley-completo-uniforme}]. Finally, we show that the result of Ishaq \cite[Theorem~2.8]{Ishaq} follows from the Theorem \ref{teorema2}.


\section{Algebraic and combinatorial Stanley depth}\label{seccion2}

For a positive integer $n$, let $[n] = \{1,...,n\}$ and let \textbf{$2^{n}$} denote the Boolean algebra consisting of all subsets of $[n]$. For $x\leq y$ in a poset $P$, we let $[x,y]$ $=$ $\{z : x\leq z\leq y\}$ and call
$[x,y]$ an interval in $P$. If $P$ is a poset and $x\in P$, we let $U[x] = \{y\in P : y\geq x\}$ and call this the up-set of $x$. In \cite{Herzog}, Herzog et al. introduced a powerful connection between the Stanley depth of a monomial ideal and a combinatorial partitioning problem for partially ordered sets. For $c\in\mathbb{N}^{n}$, let $x^{c}:=x_{1}^{c(1)}x_{2}^{c(2)}\cdots x_{n}^{c(n)}$. Let $I=(x^{v_{1}},\ldots,x^{v_{q}})\subset R$ be a monomial ideal. Let $h\in\mathbb{N}^{n}$ be such that $h\geq v_{i}$ for all $i$. The characteristic poset of $I$ with respect to $h$, denoted $P_{I}^{h}$ is defined as the induced subposet of $\mathbb{N}^{n}$ with ground set
$$\{c\in\mathbb{N}^{n} : c\leq h \textmd{  and there is } i\textmd{  such that  }c\geq v_{i}\}.$$

Let $\mathcal{D}$ be a partition of $P_{I}^{h}$ into intervals. For $J = [x,y]\in \mathcal{D}$, define
$$Z_{J}=\{i\in[n] : y(i)=h(i)\}.$$
Define the Stanley depth of a partition $\mathcal{D}$ to be
$$sdepth(\mathcal{D}) = min_{J\in\mathcal{D}}|Z_{J}|$$
and the Stanley depth of the poset $P_{I}^{h}$ to be
$$sdepth(P_{I}^{h}) = max_{\mathcal{D}}\,\,sdepth(\mathcal{D}),$$
where the maximum is taken over all partitions $\mathcal{D}$ of $P_{I}^{h}$ into intervals. Herzog et al. showed in \cite{Herzog} that
\begin{equation}\label{ideal-partition}
sdepth(I) = sdepth(P_{I}^{h}).
\end{equation}

If $I$ is a squarefree monomial ideal, then we may take $h = (1,1,\ldots ,1)$ and work inside $\{0,1\}^{n}$, which is isomorphic to \textbf{$2^{n}$}. A monomial $m$ in $R$ then can be identified with the subset of $[n]$ whose elements
correspond to the subscripts of the variables appearing in $m$. Let $G(I) = \{x^{v_{1}},x^{v_{2}},\ldots,x^{v_{q}}\}$ be the
set of minimal monomial generators of $I$ and $A_{i}\subseteq[n]$ corresponds to $v_{i}$.
The characteristic poset of $I$ with respect
to $h=(1,1,\ldots ,1)$, denoted by $P_{I}^{h}$ is in fact the set
$$P_{I}^{h}=\{C\subset[n] : C \textmd{ contains the } supp(v_{i}) \textmd{ for some } i\}$$
where $supp(v_{i}) = \{j : x_{j} | u_{i}\}$.
Then the definition of $P_{I}^{h}$ clearly simplifies to
$$P_{I}^{h}=\cup_{i=1}^{q}U[A_{i}]$$
as a subposet of \textbf{$2^{n}$}. For an interval $J = [X,Y]$, we then have that $|Z_{J}|$ corres\-ponds to $|Y|$.
\\

Let $\mathcal{P} : P_{I}^{h} = \cup_{i=1}^{q}[C_{i},D_{i}]$ be a partition of $P_{I}^{h}$, and for each $i$, let $c_{i}\in \{0,1\}^{n}$ be the $n$--tuple
such that $supp(x^{c_{i}})=C_{i}$. Then there is a Stanley decomposition $\mathcal{D(P)}$ of $I$
$$\mathcal{D(P)} : I=\bigoplus_{i=1}^{q}x^{c_{i}}K[\{x_{k} : k\in D_{i}\}].$$

The above description of $sdepth(I) = sdepth(P_{I}^{h})$ shows that
\begin{Lemma}\label{lemita}
If $I$ is a squarefree monomial ideal and $G(I)$ is the minimal monomial
generating set of $I$, then $\textmd{min}\{deg(v) : v\in G(I)\}\leq sdepth(I)\leq n$.
\end{Lemma}
By the previous lemma, if $\mathcal{C}$ is a $d$--uniform clutter, then $d\leq sdepth(I(\mathcal{C}))\leq n$.

\section{Stanley depth of edge ideals}\label{seccion3}
Let $R=K[x_1,\ldots,x_n]$ be a polynomial ring over a field $K$ and
let $v_1,\ldots,v_q$ be the column vectors of a matrix $A=(a_{ij})$ whose entries are
non-negative integers. For technical reasons, we shall always assume that the rows and
columns of the matrix $A$ are different from zero. As usual we
use the notation $x^a:=x_1^{a_1} \cdots x_n^{a_n}$,
where $a=(a_1,\ldots,a_n)\in \mathbb{N}^n$.
\\
Consider the {\it monomial ideal}:
$$
I=(x^{v_1},\ldots,x^{v_q})\subset R
$$
generated by $F=\{x^{v_1},\ldots,x^{v_q}\}$.
\\

Let $A$ be the {\it incidence matrix\/} of $\mathcal{C}$ whose column
vectors are $v_1,\ldots,v_q$. The {\it set covering polyhedron\/}
of $\mathcal{C}$ is given by:
$$
Q(A)=\{x\in\mathbb{R}^n\vert\, x\geq 0;\, xA\geq{\mathbf 1}\},
$$
A subset $C\subset V(\mathcal{C})$ is called a
{\it minimal vertex cover\/} of $\mathcal{C}$ if:
(i) every edge of $\mathcal{C}$ contains at least one vertex of $C$,
and (ii) there is no proper subset of $C$ with the first
property. The map $C\mapsto \sum_{x_i\in C}e_i$ gives a bijection
between the mi\-ni\-mal vertex covers of $\mathcal{C}$ and the integral
vectors of $Q(A)$ \cite{tesis-docto-dupont}. A polyhedron is called an {\it
integral polyhedron\/} if it has only integral vertices.

\begin{Definition}{\rm A $d$--uniform clutter $\mathcal{C}(V,E)$ with vertex set $V$ and edge set $E$ is called {\it
$k$--partite\/} if the vertex set $V$ is partitioned into $k$ disjoint subset $V_{1}, V_{2},... , V_{k}$ and
$\mid e\cap V_{i}\mid\leq1$ for all $e\in E$ and $1\leq i\leq k$.}
\end{Definition}

\begin{Definition}{\rm A $d$--uniform clutter $\mathcal{C}(V,E)$ with vertex set $V$ and edge set $E$ is called {\it complete $k$--partite\/} ($d\leq k\leq n$) if the vertex set $V$ is partitioned into $k$ disjoint subset $V_{1}, V_{2},\ldots , V_{k}$ and
$E=\{\{x_{j_{1}},\ldots ,x_{j_{d}}\} : |x_{j_{l}}\cap V_{i}|\leq1 \}$, in that case we say that $V=V_{1}\cup V_{2}\cup\cdots\cup V_{k}$ is a {\it complete partition\/}. Note that if $d=k$, then $V_{1}, V_{2},\ldots , V_{k}$ are the minimal vertex covers of $\mathcal{C}$.}
\end{Definition}

Let $I\subset R$ be the edge ideal of a complete bipartite graph over $n$ vertices with
$n\geq4$. In \cite{Ishaq} Ishaq showed that $$sdepth(I)\leq\frac{n+2}{2}.$$

Now let $\mathcal{C}$ be a complete $k$--partite $d$--uniform clutter with vertex set $V(\mathcal{C})$ partitioned into $k$ disjoint subset $V_{1}, V_{2},\ldots , V_{k}$; $V(\mathcal{C})=V_{1}\cup\cdots\cup V_{k}$, with $|V_{i}|=r_{i}$, where $r_{i}\in\mathbb{N}$ and $2\leq r_{1}\leq\cdots\leq r_{k}$. Let $r_{1}+\cdots+r_{k}=n$ and $V_{1}=\{x_{1},\ldots,x_{r_{1}}\}$, $V_{2}=\{x_{r_{1}+1},\ldots,x_{r_{1}+r_{2}}\},$ $\ldots,$
$V_{k}=\{x_{r_{1}+\cdots r_{k-1}+1},\ldots,x_{r_{1}+\cdots r_{k}}\}$. Let $I_{1}=(V_{1})$, $I_{2}=(V_{2}),$ $\ldots,$ $I_{k}=(V_{k})$ be the monomial ideals in $R$. Note that

$$|E(\mathcal{C})|=\left(\sum_{1\leq j_{1}<j_{2}<\cdots<j_{d}\leq k}r_{j_{1}}r_{j_{2}}\cdots r_{j_{d}}\right).$$
Then the edge ideal of $\mathcal{C}$ is of the form
$$ I=\sum_{1\leq j_{1}<j_{2}<\cdots<j_{d}\leq k}I_{j_{1}}\cap I_{j_{2}}\cap\cdots\cap I_{j_{d}}.$$

The next result follows from the fact that the Stanley's conjecture holds for
finite products of monomial prime ideals (see from \cite[Corollary~2.9]{pournaki}); for convenience we include
a short proof.

\begin{Theorem}\label{conjetura-stanley-completo-uniforme}
Let $I$ be the edge ideal of $d$--uniform complete $d$--partite clutter. Then Stanley's
Conjecture holds for $I$.
\end{Theorem}
\demo
We continue to use the notation used in the above description of $I=I(\mathcal{C})=(V_{1})\cap\cdots\cap(V_{d})$, with $V(\mathcal{C})=V_{1}\cup\cdots\cup V_{d}$. In our situation $V_{1}, \ldots, V_{d}$ are the minimal vertex covers of $\mathcal{C}$. Therefore $$I=(V_{1})\cap\cdots\cap(V_{d})$$
is a reduced intersection of monomial prime ideals of $R$, where $(V_{i})\nsubseteq\sum_{j=1,j\neq i}^{d}(V_{j})$ for all $1\leq i\leq d$. Then by \cite[Theorem~3.3]{Popescu1},
$$depth(I)=d\leq sdepth(I).$$
\QED

\begin{Theorem}\label{teorema1}
Let $\mathcal{C}$ be a $d$--uniform complete $k$--partite clutter. Then
$$d\leq sdepth(I(\mathcal{C}))\leq d + \frac{1}{|E(\mathcal{C})|}\left(\sum_{1\leq j_{1}<j_{2}<\cdots<j_{d}\leq k}\left(\sum_{i=1}^{d}(\begin{array}{c}
                                                                                    r_{j_{i}} \\
                                                                                    2
                                                                                  \end{array})\frac{r_{j_{1}}\cdots r_{j_{d}}}{r_{j_{i}}}
\right)\right).$$
\end{Theorem}
\demo
Note that $I=I(\mathcal{C})$ is a squarefree monomial ideal generated by monomials of degree $d$.
Let $\rho = sdepth(I)$ and $\mathcal{P} : P_{I}^{h} = \cup_{i=1}^{q}[C_{i},D_{i}]$ be a partition of $P_{I}^{h}$
satisfying $sdepth(\mathcal{D(P)})=\rho$, where $\mathcal{D(P)}$ is the Stanley decomposition of $I$ with
respect to the partition $\mathcal{P}$. We may choose $\mathcal{P^{'}}$ such that $|D|=\rho$ whe\-ne\-ver $C\neq D$ in the interval $[C,D]$,
considering these intervals of $\mathcal{P}$ with $|D|=\rho$ and $1$-dimensional spaces. Now we see that for each interval $[C,D]$ in $\mathcal{P^{'}}$ with $|C|=d$
we have $\rho - d$ subsets of cardinality $d+1$ in this interval. The total number of these
kind of intervals is $|E(\mathcal{C})|=\sum r_{j_{1}}r_{j_{2}}\cdots r_{j_{d}}$, where the sum runs over all $1\leq j_{1}<j_{2}<\cdots<j_{d}\leq k$. So we have $$(\rho-d)\left(\sum r_{j_{1}}r_{j_{2}}\cdots r_{j_{d}}\right)$$
subsets of cardinality $d+1$. This number is less than or equal to the total number of
monomials $m\in I$ with $\deg (m)=d+1$ and $|supp(m)|=d+1$. Furthermore,
$$\{m : m\in I; \deg (m)=d+1; supp(m)=d+1\}=\{x^{e}x_{i} : e\in E(\mathcal{C}); i\notin e\},$$
with cardinality $\sum_{1\leq j_{1}<j_{2}<\cdots<j_{d}\leq k}\left(\sum_{i=1}^{d}(\begin{array}{c}
                                                                                    r_{j_{i}} \\
                                                                                    2
                                                                                  \end{array})\frac{r_{j_{1}}\cdots r_{j_{d}}}{r_{j_{i}}}
\right)$.
Hence
$$(\rho-d)|E(\mathcal{C})| \leq \sum_{1\leq j_{1}<j_{2}<\cdots<j_{d}\leq k}\left(\sum_{i=1}^{d}(\begin{array}{c}
                                                                                    r_{j_{i}} \\
                                                                                    2
                                                                                  \end{array})\frac{r_{j_{1}}\cdots r_{j_{d}}}{r_{j_{i}}}
\right).$$
Therefore we obtain the required inequality.
\QED
\begin{Theorem}\label{teorema1.seguido}
Let $\mathcal{C}$ be a $d$--uniform complete $d$--partite clutter. Then
$$d\leq sdepth(I(\mathcal{C}))\leq d + \sum_{i=1}^{d}\frac{r_{i}-1}{2}.$$
\end{Theorem}
\demo
The proof is analogous to the proof of Theorem~\ref{teorema1}, but with $|E(\mathcal{C})|=r_{1}r_{2}\cdots r_{d}$ and
$$\{m : m\in I; \deg (m)=d+1; supp(m)=d+1\}=\{x^{e}x_{i} : e\in E(\mathcal{C}); i\notin e\},$$
has cardinality $$\sum_{i=1}^{d}\left(\begin{array}{c}
                                        r_{i} \\
                                        2
                                      \end{array}\right)
\frac{r_{1}r_{2}\cdots r_{d}}{r_{i}} = (r_{1}r_{2}\cdots r_{d})\sum_{i=1}^{d}\frac{r_{i}-1}{2}.$$
Hence
$$(\rho-d)(r_{1}r_{2}\cdots r_{d})\leq (r_{1}r_{2}\cdots r_{d})\sum_{i=1}^{d}\frac{r_{i}-1}{2}.$$
Therefore we obtain the required inequality.
\QED

\begin{Definition}{\rm A clutter $\mathcal{C}(V,E)$, whose set covering polyhedron
$Q(A)$ is integral, is called {\it
integral\/}.}
\end{Definition}
\begin{Lemma}\label{lemma-1}\emph{(See \cite{tesis-docto-dupont})}
If $\mathcal{C}$ is an integral $d$--uniform clutter,
then there exists a mi\-ni\-mal vertex cover of $\mathcal{C}$ intersecting every edge of $\mathcal{C}$ in
exactly one vertex.
\end{Lemma}
\demo
Let $B$ be the integral matrix whose columns are the
vertices of $Q(A)$. It is not hard to show that {\it a vector
$\alpha\in\mathbb{R}^n$ is an integral vertex of $Q(A)$ if and only if
$\alpha=\sum_{x_i\in C}e_i$ for some minimal vertex cover
$C$ of $\mathcal{C}$}. Thus the columns of $B$ are the characteristic
vectors of the minimal vertex covers of $\mathcal C$.
Using \cite[Theorem~1.17]{cornu-book} we get that
$$
Q(B)=\{x\vert\,
x\geq 0;xB\geq\mathbf{1}\}
$$
is an integral
polyhedron whose vertices are the columns of $A$, where $\mathbf{1}=(1,1,\ldots,1)$. Therefore we have
the equality
\begin{equation}\label{may10-09}
Q(B)=\mathbb{R}_+^n+{\rm conv}(v_1,\ldots,v_q).
\end{equation}

We proceed by contradiction. Assume that for each column $u_k$
of $B$ there
exists a vector $v_{i_k}$ in $\{v_1,\ldots,v_q\}$ such that
$\langle v_{i_k},u_k\rangle\geq 2$. Here $\langle\,   ,  \rangle$ is the standard inner product. Then
$$
v_{i_k}B\geq \mathbf{1}+e_k,
$$ where $e_{i}$ is the $i$-th unit vector.
\\
Consider the vector $\alpha=v_{i_1}+\cdots+v_{i_s}$, where $s$ is the
number of columns of $B$. From the
inequality
$$
{\alpha}B\geq (\mathbf{1}+e_1)+\cdots+(\mathbf{1}+e_s)=(s+1,\ldots,s+1)
$$
we obtain that $\alpha/(s+1)\in Q(B)$. Thus, using
Eq.~(\ref{may10-09}),
we can write
\begin{equation}\label{may19-09}
\alpha/(s+1)=\mu_1e_1+\cdots+\mu_ne_n+\lambda_1v_1+\cdots+\lambda_qv_q\
\ \ \ (\mu_i,\lambda_j\geq 0;\ \textstyle\sum\lambda_i=1).
\end{equation}
Therefore taking inner products with $\mathbf{1}$ in
Eq.~(\ref{may19-09}) and using the fact that $\mathcal{C}$ is $d$--uniform
we get that $|\alpha|\geq (s+1)d$. Then using the equality
$\alpha=v_{i_1}+\cdots+v_{i_s}$ we conclude
$$
sd=|v_{i_1}|+\cdots+|v_{i_s}|=|\alpha|\geq (s+1)d,
$$
a contradiction because $d\geq 1$.
\QED
\\
\\

Let $\mathcal{C}$ be a clutter and let $I=I(\mathcal{C})$ be its edge
ideal. Recall that a {\it deletion\/} of
$I$ is any ideal $I'$ obtained from $I$
by making a variable equal to $0$. A {\it deletion\/} of $\cal C$ is
a clutter ${\cal C}'$ that
corresponds to a deletion $I'$ of $I$.  Notice that
${\cal C}'$ is
obtained from $I'$ by considering the unique set
of square-free monomials that minimally generate $I'$. A {\it
contraction\/} of
$I$ is any ideal $I'$ obtained from $I$
by making a variable equal to $1$. A {\it contraction\/} of $\cal C$ is
a clutter ${\cal C}'$ that corresponds to a contraction $I'$ of $I$. This terminology is
consistent with that of \cite[p.~23]{cornu-book}.
\\
A clutter obtained from $\mathcal{C}$ by a sequence of deletions and
contractions of vertices is called a {\it minor\/} of $\mathcal{C}$. The clutter $\mathcal{C}$ is considered itself a minor.
\\
\\
The notion of a minor of a clutter is not a generalization of the notion of a
minor of a graph in the sense of graph theory \cite[p.~25]{Schr2}.
For instance if $G$ is a cycle of length four and we contract an edge we obtain that
a triangle is a minor of $G$, but a triangle cannot be a minor of $G$
in our sense.
\\

The notion of a minor plays a prominent role in
combinatorial optimization \cite{cornu-book}. As an application of the power of using minors, this allows us to get a nice decomposition of an integral uniform clutter.
\begin{Proposition}\label{propos}\emph{(See \cite{tesis-docto-dupont})}
If $\mathcal{C}(V,E)$ be an integral $d$--uniform clutter, then there are $V_{1}, ..., V_{d}$ mutually disjoint minimal vertex covers of $\mathcal{C}$ such that $V = \bigcup_{i=1}^{d}V_{i}$. In particular $|supp(x^{e})\cap V_{i}| = 1$ for all
$e\in E; 1\leq i\leq d$.
\end{Proposition}
\demo
 By induction on $d$. If $d=1$, then
$E(\mathcal{C})=\{\{x_1\},\ldots,\{x_n\}\}$ and $V$ is a minimal
vertex cover of $\mathcal{C}$. In this case we set $V_1=V$.
Assume $d\geq 2$. By
Lemma~\ref{lemma-1} there is a minimal vertex
cover $V_1$ of $\mathcal C$ such that $|{\rm supp}(x^{v_i})\cap
V_1|=1$ for all
$i$. Consider the ideal $I'$ obtained from $I$ by making $x_i=1$ for
$x_i\in V_1$. Let $\mathcal{C}'$ be the clutter corresponding to $I'$
and let $A'$ be the incidence matrix of $\mathcal{C}'$. The ideal $I'$
(resp. the clutter $\mathcal{C}'$) is a minor of $I$ (resp.
$\mathcal{C}$). Recall that the integrality of $Q(A)$ is preserved
under taking
minors \cite[Theorem~78.2]{Schr2}, so $Q(A')$ is integral. Then
$\mathcal{C}'$
is a $(d-1)$-uniform clutter whose set covering polyhedron $Q(A')$ is
integral. Note that $V(\mathcal{C}')=V\setminus V_1$.
Therefore by induction hypothesis there are $V_2,\ldots,V_d$
pairwise disjoint minimal vertex covers of $\mathcal{C}'$ such that
$V\setminus V_1=V_2\cup\cdots\cup V_d$. To complete the proof
observe that $V_2,\ldots,V_d$ are minimal vertex covers of
$\mathcal{C}$. Indeed if $e$ is an edge of $\mathcal{C}$ and $2\leq
k\leq d$, then $e\cap V_1=\{x_i\}$ for some $i$. Since
$e\setminus\{x_i\}$ is an edge of $\mathcal{C}'$, we get
$(e\setminus\{x_i\})\cap V_k\neq\emptyset$. Hence $V_k$ is a vertex cover of
$\mathcal{C}$. Furthermore if $x\in V_k$, then by the minimality of
$V_k$ relative to $\mathcal{C}'$ there is an edge $e'$ of
$\mathcal{C}'$ disjoint from
$V_k\setminus\{x\}$.
Since $e=e'\cup \{y\}$ is an edge of $\mathcal{C}$ for some $y\in
V_1$, we obtain that $e$ is an edge of $\mathcal{C}$ disjoint from
$V_k\setminus\{x\}$. Therefore $V_k$ is a minimal vertex cover of
$\mathcal{C}$, as required.
\QED
\begin{Theorem}\label{teorema2}
Let $\mathcal{C}(V,E)$ be an integral $d$--uniform clutter. Then $\mathcal{C}$ is a $d$--partite clutter, with
$$d\leq sdepth(I(\mathcal{C}))\leq d + \frac{r_{1}\cdots r_{d}}{|E(\mathcal{C})|}\sum_{i=1}^{d}\frac{r_{i}-1}{2}.$$
\end{Theorem}
\demo
By Proposition~\ref{propos}, we have that $\mathcal{C}$ is a $d$--partite clutter. The proof is analogous to the proof of Theorem~\ref{teorema1}, but with $k=d$.
Note that $I=I(\mathcal{C})$ is a squarefree monomial ideal generated by monomials of degree $d$.
Let $\rho = sdepth(I)$ and $\mathcal{P} : P_{I}^{h} = \cup_{i=1}^{q}[C_{i},D_{i}]$ be a partition of $P_{I}^{h}$
satisfying $sdepth(\mathcal{D(P)})=\rho$, where $\mathcal{D(P)}$ is the Stanley decomposition of $I$ with
respect to the partition $\mathcal{P}$. We may choose $\mathcal{P}$ such that $|D|=\rho$ whenever $C\neq D$
in the interval $[C,D]$. Now we see that for each interval $[C,D]$ in $\mathcal{P}$ with $|C|=d$
we have $\rho - d$ subsets of cardinality $d+1$ in this interval. The total number of these
kind of intervals is $|E(\mathcal{C})|$. 
So we have $$(\rho-d)|E(\mathcal{C})|$$  
subsets of cardinality $d+1$. This number is less than or equal to the total number of
monomials $m\in I$ with $\deg (m)=d+1$ and $|supp(m)|=d+1$. Furthermore,
$$\{m : m\in I; \deg (m)=d+1; supp(m)=d+1\}=\{x^{e}x_{i} : e\in E(\mathcal{C}); i\notin e\},$$
with cardinality less than or equal $\sum_{i=1}^{d}(\begin{array}{c}
                                                                                    r_{i} \\
                                                                                    2
                                                                                  \end{array})\frac{r_{1}\cdots r_{d}}{r_{i}}
                                                                                  = \sum_{i=1}^{d}\frac{r_{i}-1}{2}r_{1}\cdots r_{d}$.
Therefore we obtain
$$sdepth(I(\mathcal{C}))\leq d + \frac{1}{|E(\mathcal{C})|}\left(\sum_{i=1}^{d}\left(\frac{r_{i}-1}{2}\right)r_{1}\cdots r_{d}\right).$$
Hence
$$sdepth(I(\mathcal{C}))\leq d + \frac{r_{1}\cdots r_{d}}{|E(\mathcal{C})|}\sum_{i=1}^{d}\frac{r_{i}-1}{2}.$$
\QED
\begin{Corollary}
Let $\mathcal{C}(V,E)$ be an integral $d$--uniform clutter, such that its decomposition $d$--partite $V=V_{1}\cup V_{2}\cup\cdots\cup V_{k}$ is complete. Then
$$d\leq sdepth(I(\mathcal{C}))\leq d + \sum_{i=1}^{d}\frac{r_{i}-1}{2}.$$
\end{Corollary}
\demo
It follows from Theorem~\ref{teorema1.seguido} or Theorem~\ref{teorema2}.
\QED

\begin{Corollary}\emph{(\cite[Theorem~2.8]{Ishaq})}
The Stanley depth of the edge ideal of a complete bipartite graph over $n$ vertices with $n\geq4$ is less than or equal to $\frac{n+2}{2}$.
\end{Corollary}
\demo
This follows from the fact that complete bipartite graphs are integral clutters.
\QED

\bibliographystyle{plain}

\begin{thebibliography}{10}

\bibitem{Apel} Apel, J., On a conjecture of R.P. Stanley, Part I — Monomial ideals, J. Algebraic Combin.
(2003) {\bf 17},  39--56.

\bibitem{Anwar-Popescu} Anwar, I. and Popescu D., Stanley conjecture in small embedding dimension, J. Algebra (2007) {\bf 318}, 1027--1031.

\bibitem{Biro} Bir\'o, C., Howard, D. M., Keller, M. T., Trotter,  W. T., Young,  S. J., Interval partitions and Stanley depth, J. Combin. Theory Ser. A  (2010) {\bf 117},  no. 4, 475–-482.


\bibitem{Cimpoeas} Cimpoea\c{s}, M., Stanley depth of complete intersection monomial ideals, Bull. Math. Soc. Sci. Math. Roumanie (N.S.)
(2008) {\bf 51 (99)} (3), 205--211.

\bibitem{cornu-book} {Cornu\'ejols, G. {\it Combinatorial optimization:
Packing and covering\/}, CBMS-NSF Regional Conference Series in Applied
Mathematics (2001) {\bf 74}, SIAM.}


\bibitem{tesis-docto-dupont} {Dupont, L. A., {\it Rees Algebras, Monomial Subrings and Linear Optimization Problems\/}, PhD Thesis (2010), Cinvestav-IPN. arXiv:1006.2774.}

\bibitem{normali} Escobar, C., Villarreal,  R. H., Yoshino, Y., Torsion
freeness and normality of blowup rings of monomial ideals,
{\it Commutative Algebra\/}, Lect. Notes Pure Appl. Math. (2006)
{\bf 244}, Chapman \& Hall/CRC, Boca Raton, FL, pp. 69-84.

\bibitem {faridi} Faridi, S., The facet ideal of a simplicial complex,
Manuscripta Math. (2002) {\bf 109}, 159-174.

\bibitem{Herzog} Herzog, J., Vl¢adoiu, M., Zheng, X., How to compute the Stanley depth of a monomial ideal, J.
Algebra (2009) {\bf 322}, 3151--3169.

\bibitem{Ishaq} Ishaq, M., Upper bounds for the Stanley depth, Comm. Algebra  (2012) {\bf 40},  no. 1, 87–-97.

\bibitem{Ishaq2} Ishaq, M., Qureshi, M. I., Stanley depth of edge ideals, Studia Sci. Math. Hungar.  (2012) {\bf 49},  no. 4, 501–-508.

\bibitem{Popescu2} Popescu, D., Stanley depth of multigraded modules, J. Algebra (2009) {\bf 321}, 2782--2797.

\bibitem{Popescu3} Popescu, D., Qureshi, M. I., Computing the Stanley depth, J. Algebra (2010) {\bf 323}, 2943--2959.

\bibitem{Popescu1} Popescu, D., Stanley conjecture on intersection of four monomial prime ideals,
Comm. Algebra  (2013) {\bf 41},  no. 11, 4351–-4362.

\bibitem{pournaki} Pournaki M. R., Seyed Fakhari S. A. and Yassemi S., On the Stanley depth of weakly polymatroidal ideals, Arch. Math. (2013) {\bf 100}, no. 2, 115--121.

\bibitem{Schr2} {Schrijver, A., {\it Combinatorial Optimization\/},
Algorithms and Combinatorics (2003) {\bf 24}, Springer-Verlag, Berlin.}

\bibitem{Shen} Shen, Y. H., Stanley depth of complete intersection monomial ideals and upper-discrete partitions, J. Algebra (2009) 321, (4) 1285--1292.

\bibitem{ITG}{Simis, A., Vasconcelos, W.~V., Villarreal, R. H., On
the ideal theory of graphs, J. Algebra, (1994) {\bf 167}, 389--416.}

\bibitem{conj-Stanley} Stanley, R. P., Linear Diophantine equations and local cohomology, Invent. Math. (1982) {\bf 68}, 175--193.

\bibitem{Stanley} Stanley, R. P., {\it Combinatorics and Commutative
Algebra. Second edition.}
Progress in Mathematics (1996) {\bf 41}. Birkh{\"a}user Boston, Inc.,
Boston, MA.

\bibitem{Vi2}{Villarreal, R. H., Cohen-{M}acaulay graphs, Manuscripta
Math. (1990) {\bf 66}, 277--293.}

\bibitem{monalg}{ Villarreal, R. H., {\it Monomial
Algebras\/}, Dekker, New York, N.Y.} (2001).


\bibitem{zheng} Zheng, X., Resolutions of facet ideals,
Comm. Algebra (2004) {\bf 32 (6)},  2301--2324.




\end{thebibliography}

\end{document}